\documentclass[12pt,reqno]{amsart} 
\usepackage{amssymb,amscd}

\begin{document}

%%%%%%%%%%%%% Set date and title %%%%%%%%%%%%%%%%%%%%%%
% \newif\ifdraft \drafttrue
\newcommand{\DATE}{December 2003}
\newcommand{\TITLE}{A lower bound for the canonical height on abelian varieties over abelian extensions}
\newcommand{\TITLERUNNING}{Heights on abelian varieties over abelian extensions}
%%%%%%%%%%%%%%%%%%%%%%%%%%%%%%%%%%%%%%%%%%%%%%%%%%%%%%%%%%%%%%%%%%%%%%

% Page settings
%%   \topmargin = 0 in
%%   \headsep = .1 in
%%   \textwidth = 6.5 in
%%   \textheight = 8.9 in
%%   \baselineskip = .16666 in
%%   \oddsidemargin = 0 in
%%   \evensidemargin = 0 in

%%%%%%%%%%%%%%%%%%%%%%%%%%%%%%%%%%%%%%%%%%%%%%%%%%%%%%%%%%%%%%%%%%%%%%
% Theorem environments

\numberwithin{equation}{section}

\newtheorem{theorem}[equation]{Theorem}
\newtheorem{lemma}[equation]{Lemma}
\newtheorem{conjecture}[equation]{Conjecture}
\newtheorem{proposition}[equation]{Proposition}
\newtheorem{corollary}[equation]{Corollary}

\theoremstyle{definition}
% The * surpresses numbering, so definitions are not numbered
\newtheorem*{definition}{Definition}
\newtheorem{example}[equation]{Example}

\theoremstyle{remark}
\newtheorem{remark}[equation]{Remark}
\newtheorem{remarks}[equation]{Remarks}
\newtheorem*{acknowledgement}{Acknowledgements}

%%%%%%%%%%%%%%%%%%%%%%%%%%%%%%%%%%%%%%%%%%%%%%%%%%%%%%%%%%%%%%%%%%%%%%

%%%%%%%% Set Up Environment for Notation %%%%%%%%%%%%%%
% This is currently set to allow quite wide items to be defined
\newenvironment{notation}[0]{%
  \begin{list}%
    {}%
    {\setlength{\itemindent}{0pt}
     \setlength{\labelwidth}{4\parindent}
     \setlength{\labelsep}{\parindent}
     \setlength{\leftmargin}{5\parindent}
     \setlength{\itemsep}{0pt}
     }%
   }%
  {\end{list}}

%%%%%%%% Set Up Environment for Parts in Theorems %%%%%%%%%%%%%%
\newenvironment{parts}[0]{%
  \begin{list}{}%
    {\setlength{\itemindent}{0pt}
     \setlength{\labelwidth}{1.5\parindent}
     \setlength{\labelsep}{.5\parindent}
     \setlength{\leftmargin}{2\parindent}
     \setlength{\itemsep}{0pt}
     }%
   }%
  {\end{list}}
% Use \Part{(a)}, instead of \item[(a)], to ensure upright font
\newcommand{\Part}[1]{\item[\upshape#1]}

%%%%%%%%%%%%%%%%%%
% Greek Alphabet %
%%%%%%%%%%%%%%%%%%
\renewcommand{\a}{\alpha}
\renewcommand{\b}{\beta}
\newcommand{\g}{\gamma}
\renewcommand{\d}{\delta}
\newcommand{\e}{\epsilon}
\newcommand{\f}{\phi}
\renewcommand{\l}{\lambda}
\renewcommand{\k}{\kappa}
\newcommand{\lhat}{\hat\lambda}
\newcommand{\m}{\mu}
\renewcommand{\o}{\omega}
\renewcommand{\r}{\rho}
\newcommand{\rbar}{{\bar\rho}}
\newcommand{\s}{\sigma}
\newcommand{\sbar}{{\bar\sigma}}
\renewcommand{\t}{\tau}
\newcommand{\z}{\zeta}

\newcommand{\D}{\Delta}
\newcommand{\F}{\Phi}

%%%%%%%%%%%%%%%%%%%%
% Fraktur Alphabet %
%%%%%%%%%%%%%%%%%%%%
\newcommand{\gp}{{\mathfrak{p}}}
\newcommand{\gP}{{\mathfrak{P}}}
\newcommand{\gq}{{\mathfrak{q}}}
\newcommand{\gf}{{\mathfrak{f}}}

%%%%%%%%%%%%%%%%%%%%%%%%%
% Calligraphic Alphabet %
%%%%%%%%%%%%%%%%%%%%%%%%%
\newcommand{\Acal}{{\mathcal A}}
\newcommand{\Bcal}{{\mathcal B}}
\newcommand{\Ccal}{{\mathcal C}}
\newcommand{\Dcal}{{\mathcal D}}
\newcommand{\Ecal}{{\mathcal E}}
\newcommand{\Fcal}{{\mathcal F}}
\newcommand{\Gcal}{{\mathcal G}}
\newcommand{\Hcal}{{\mathcal H}}
\newcommand{\Ical}{{\mathcal I}}
\newcommand{\Jcal}{{\mathcal J}}
\newcommand{\Kcal}{{\mathcal K}}
\newcommand{\Lcal}{{\mathcal L}}
\newcommand{\Mcal}{{\mathcal M}}
\newcommand{\Ncal}{{\mathcal N}}
\newcommand{\Ocal}{{\mathcal O}}
\newcommand{\Pcal}{{\mathcal P}}
\newcommand{\Qcal}{{\mathcal Q}}
\newcommand{\Rcal}{{\mathcal R}}
\newcommand{\Scal}{{\mathcal S}}
\newcommand{\Tcal}{{\mathcal T}}
\newcommand{\Ucal}{{\mathcal U}}
\newcommand{\Vcal}{{\mathcal V}}
\newcommand{\Wcal}{{\mathcal W}}
\newcommand{\Xcal}{{\mathcal X}}
\newcommand{\Ycal}{{\mathcal Y}}
\newcommand{\Zcal}{{\mathcal Z}}

\newcommand{\OO}{{\mathcal O}}    % Structure sheaf
\renewcommand{\O}{{\mathcal O}}   % Ring of integers

%%%%%%%%%%%%%%%%%%%%%%%%%%%%
% Blackboard Bold Alphabet %
%%%%%%%%%%%%%%%%%%%%%%%%%%%%
\renewcommand{\AA}{\mathbb{A}}
\newcommand{\BB}{\mathbb{B}}
\newcommand{\CC}{\mathbb{C}}
\newcommand{\FF}{\mathbb{F}}
\newcommand{\GG}{\mathbb{G}}
\newcommand{\PP}{\mathbb{P}}
\newcommand{\QQ}{\mathbb{Q}}
\newcommand{\RR}{\mathbb{R}}
\newcommand{\ZZ}{\mathbb{Z}}

%%%%%%%%%%%%%%%%%%%%%%%%%%
% Boldface Math Alphabet %
%%%%%%%%%%%%%%%%%%%%%%%%%%
\newcommand{\bfa}{{\mathbf a}}
\newcommand{\bfb}{{\mathbf b}}
\newcommand{\bfc}{{\mathbf c}}
\newcommand{\bfe}{{\mathbf e}}
\newcommand{\bff}{{\mathbf f}}
\newcommand{\bfg}{{\mathbf g}}
\newcommand{\bfp}{{\mathbf p}}
\newcommand{\bfr}{{\mathbf r}}
\newcommand{\bfs}{{\mathbf s}}
\newcommand{\bft}{{\mathbf t}}
\newcommand{\bfu}{{\mathbf u}}
\newcommand{\bfv}{{\mathbf v}}
\newcommand{\bfw}{{\mathbf w}}
\newcommand{\bfx}{{\mathbf x}}
\newcommand{\bfy}{{\mathbf y}}
\newcommand{\bfz}{{\mathbf z}}
\newcommand{\bfA}{{\mathbf A}}
\newcommand{\bfB}{{\mathbf B}}
\newcommand{\bfC}{{\mathbf C}}
\newcommand{\bfF}{{\mathbf F}}
\newcommand{\bfG}{{\mathbf G}}
\newcommand{\bfI}{{\mathbf I}}
\newcommand{\bfM}{{\mathbf M}}
\newcommand{\bfzero}{{\boldsymbol{0}}}
\newcommand{\bfmu}{{\boldsymbol\mu}}

%%%%%%%%%%%%%%%%%%%%%%%%%%%%%%
% Miscellaneous New Commands %
%%%%%%%%%%%%%%%%%%%%%%%%%%%%%%
\newcommand{\ab}{{\textup{ab}}}
\newcommand{\Ahat}{{\hat A}}
\newcommand{\Aut}{\operatorname{Aut}}
\newcommand{\Cond}{{\mathfrak{N}}} 
\newcommand{\Disc}{\operatorname{Disc}}
\newcommand{\Div}{\operatorname{Div}}
\newcommand{\End}{\operatorname{End}}
\newcommand{\Frob}{\operatorname{Frob}}
\newcommand{\GK}{G_{\Kbar/K}}
\newcommand{\GLK}{G_{L/K}}
\newcommand{\GL}{\operatorname{GL}}
\newcommand{\Gal}{\operatorname{Gal}}
\newcommand{\GSp}{\operatorname{GSp}}
\newcommand{\hhat}{{\hat h}}
\newcommand{\Hypothesis}{\textsc{Hyp}}
\newcommand{\Image}{\operatorname{Image}}
\newcommand{\into}{\hookrightarrow}     % an injection into
\newcommand{\Kab}{K^\ab}
\newcommand{\Kbar}{{\bar K}}
\newcommand{\Kvbar}{{{\bar K}_v}}
\newcommand{\MOD}[1]{~(\textup{mod}~#1)}
\newcommand{\Norm}{\operatorname{N}}
\newcommand{\notdivide}{\nmid}
\newcommand{\nr}{{\textup{nr}}}    % non-ramified
\newcommand{\ord}{\operatorname{ord}}
\newcommand{\Pic}{\operatorname{Pic}}
\newcommand{\Qbar}{{\bar{\QQ}}}
\newcommand{\QQbar}{{\bar{\QQ}}}
\newcommand{\rank}{\operatorname{rank}}
\newcommand{\res}{\operatornamewithlimits{res}}
\newcommand{\Resultant}{\operatorname{Resultant}}
\renewcommand{\setminus}{\smallsetminus}
\newcommand{\Spec}{\operatorname{Spec}}
\newcommand{\Support}{\operatorname{Support}}
\newcommand{\tors}{{\textup{tors}}}
\newcommand{\val}{{\operatorname{val}}}
\newcommand{\<}{\langle}
\newcommand{\la}{{\langle}}
\renewcommand{\>}{\rangle}
\newcommand{\ra}{{\rangle}}

%%%%%%%%%%%%%%%%%%%%%%%%%%%%%%%%%%%%%%%%%%%%%%%%%%%%%%%%%%%%%%%%%%%%%%

\def\BigStrut{\vphantom{$(^{(^(}_{(}$}} % Add space in tables

%%%%%%%%%%%%%%%%%%%%%%%%%%%%%%%%%%%%%%%%%%%%%%%%%%%%%%%%%%%%%%%%%%%%%%

\hyphenation{archi-me-dean}

%%%%%%%%%%%%%%%%%%%%%%%%%%%%%%%%%%%%%%%%%%%%%%%%%%%%%%%%%%%%%%%%%%%%%%

% \title[\TITLERUNNING]{\TITLE\ifdraft---Draft \DRAFTNUMBER\fi}
\title[Heights on abelian varieties over abelian extensions]
{A lower bound for the canonical height on abelian varieties over abelian extensions}
\date{\DATE}

\author[Matthew Baker]{Matthew H. Baker}
\author[Joseph Silverman]{Joseph H. Silverman}
\email{mbaker@math.uga.edu \\
       jhs@math.brown.edu}
\address{Department of Mathematics,
         University of Georgia, Athens, GA 30602-7403, USA}
\address{Mathematics Department, Box 1917
         Brown University, Providence, RI 02912 USA}

\thanks{The first author's research was supported in part
by an NSF Postdoctoral Research Fellowship, and by NSF Research Grant DMS-0300784. The second author's
research was partially supported by NSA Research Grant H98230-04-1-0064.}

\begin{abstract}
Let~$A$ be an abelian variety defined over a number field~$K$ and
let~$\hhat$ be the canonical height function on~$A(\Kbar)$ attached to a
symmetric ample line bundle~$\Lcal$. We prove that there exists a
constant~$C = C(A, K,\Lcal) > 0$ such that~$\hhat(P) \geq C$ for all
nontorsion points~$P \in A(K^\ab)$, where~$K^\ab$ is the maximal
abelian extension of~$K$.
\end{abstract}

\maketitle

%%%%%%%%%%%%%%%%%%%%%%%%%%%%%%%%%%%%%%%%%%%%%%%%%%%%%%%%%%%%%%%%%%%%%%

\section*{Introduction}
The purpose of this paper is to provide a proof of the following result.

\begin{theorem}
\label{theorem:maintheorem}
Let~$A/K$ be an abelian variety defined over a number field~$K$.  Let
$\Lcal$ be a symmetric ample line bundle on~$A/K$, and let~$\hhat :
A(\Kbar) \to \RR$ be the associated canonical height function.  Then
there exists a constant~$C = C(A, K ,\Lcal) > 0$ such that
\[
  \hhat(P) \geq C 
  \quad\text{for all nontorsion points~$P \in A(K^\ab)$.}
\]
\end{theorem}

Together with Remark~\ref{remark:Ruppert}, Theorem~\ref{theorem:maintheorem}
establishes both parts of Conjecture 1.10 in \cite{Baker}.

\par

In the case of elliptic curves, Theorem~\ref{theorem:maintheorem} was
proven by the first author~\cite{Baker} for CM elliptic curves and by
the second author~\cite{SilvermanECAbExt} for non-CM elliptic curves.
The present paper is thus a combination and extension of these two
earlier works to abelian varieties of arbitrary dimension.  The two 
papers \cite{Baker} and \cite{SilvermanECAbExt} were themselves
motivated by work of Amoroso-Dvornicich and Amoroso-Zannier
(\cite{AmorosoDvornicich}, \cite{AmorosoZannier}) providing absolute lower bounds for the height of
nontorsion points in the multiplicative group over abelian extensions of
a given number field $K$.  For more history and background
information on the elliptic curve and multiplicative group analogues
of Theorem~\ref{theorem:maintheorem}, we refer the reader to sections 1 and 2
of \cite{SilvermanECAbExt}.

\par

Theorem~\ref{theorem:maintheorem} fits into the general context of
providing lower bounds for the heights of nontorsion points on abelian
varieties.  An open problem in this context is the generalized
Lehmer conjecture (Conjecture 1.4 of \cite{DavidHindry}):

\begin{conjecture}
Suppose $A$ is an abelian variety over a number field $K$, and that
$\Lcal$ is a symmetric ample line bundle on $A$.  Then there exists 
a positive constant $c = c(A, K, \Lcal)$ such that
\[
\hhat(P) \geq c D^{-1/g_0(P)}
\]
for all nontorsion points $P \in A(\Kbar)$, where $D = [K(P):K]$ and
$g_0(P)$ denotes the dimension of the smallest algebraic subgroup of
$A$ containing $P$.
\end{conjecture}

Table~\ref{table:lowerbdinKbarhistory} gives an abbreviated
history of lower bounds for~$\hhat(P)$, where $A$ is an abelian
variety of dimension $g$ over a number field $K$, $\Lcal$ is a symmetric ample line
bundle on $A$, $P \in A(\Kbar)$ is a nontorsion point, and $D =
[K(P):K]$.  In the table, $c$ denotes a positive
constant that depends on~$A/K$ and on $\Lcal$, but not on~$P$, and
$\kappa$ denotes a positive constant depending only on the dimension
of $A$.

%%%%%%%%%%%%%%%%%%%%%%%%%%%%%%%%%%%%%%%%%%%
\begin{table}
\begin{center}
\begin{tabular}{|c|c|c|} \hline
\BigStrut
   Lower bound for $\hhat(P)$ & Restrictions & Reference \\ \hline\hline
\BigStrut
   $cD^{-\kappa}$ & none & Masser (1984) \cite{MasserAV} \\ \hline
\BigStrut
   $cD^{-\frac{1}{g}}\left(\frac{\log\log(3D)}{\log(3D)}\right)^\kappa$
     & CM, $g=g_0(P)$ & David-Hindry (2000) \cite{DavidHindry} \\ \hline
\BigStrut
   $c$ &
    $\genfrac{}{}{0pt}{0}{\text{$P \in A(\Kab)$}}{\text{$K$ totally real}}$
    & Zhang (1998)
        \cite{Zhang} \\ \hline
\BigStrut
   $c$ &
      $\genfrac{}{}{0pt}{0}{\text{$P \in A(\Kab)$}}{\text{$g=1$}}$ &
    $\genfrac{}{}{0pt}{0}{\text{Baker, Silverman}}{\text{(2003) 
\cite{Baker},\cite{SilvermanECAbExt}}}$
    \\ \hline
\end{tabular}
\caption{History of lower bounds for $\hhat$ in $A(\Kbar)$}
\label{table:lowerbdinKbarhistory}
\end{center}
\end{table}
%%%%%%%%%%%%%%%%%%%%%%%%%%%%%%%%%%%%%%%%%%%

In light of the results of \cite{AmorosoDvornicich} and \cite{AmorosoZannier} for the
multiplicative group, one might also ask if Theorem~\ref{theorem:maintheorem} admits a
generalization to arbitrary semiabelian varieties.  We do not address
this question here.

\begin{acknowledgement}
The authors would like to thank Sinnou David for suggesting
Lemma~\ref{lemma:heightvsdistance} as a means of avoiding 
the use of N{\'e}ron local heights, and Everett Howe and Mike Rosen
for their assistance with the proof of
Theorem~\ref{theorem:ordinaryCMtheorem}.
They also thank the anonymous referee for pointing out a mistake in
the original proof of Proposition 7.3.
\end{acknowledgement}

%%%%%%%%%%%%%%%%%%%%%%%%%%%%%%%%%%%%%%%%%%%%%%%%%%%%%%%%%%%%%%%%%%%%%%
\section{Notation}
\label{section:Notation}

We set the following notation and normalizations, which will be used
throughout this paper.

\begin{notation}
\item[$K$]
a number field.
\item[$K^\ab$]
the maximal abelian extension of $K$.
\item[$\O_K$]
the ring of integers of $K$.
\item[$\F_\gp$]
the residue field $\O_K/\gp$ for a prime ideal~$\gp$ of~$K$.
\item[$q_\gp$]
the absolute norm $\Norm_{K/\QQ}(\gp)$ of the ideal~$\gp$.
\item[$M_K$]
the set of places of~$K$. 
\item[$n_v$]
the normalized local degree $[K_v:\QQ_v]/[K:\QQ]$, for $v\in M_K$.
\item[$|\,\cdot\,|_v$] 
the unique absolute value on~$K$ in the equivalence class of~$v\in
M_K$ that extends the standard absolute value on the
completion~$\QQ_v$.
\item[$h$]
The absolute logarithmic Weil height
$h:\PP^n(\QQbar)\to\RR$ defined for $[x_0,\dots,x_n]\in\PP^n(K)$ by
\[
  h\bigl([x_0,\dots,x_n]\bigr) = \sum_{v \in M_K} n_v 
  \log \max \{ |x_0|_v, \ldots, |x_n|_v \}.
\]
\end{notation}

The height~$h(P)$ is well-defined, independent of the choice of a
particular number field~$K$ over which the point~$P=[x_0,\dots,x_n]$ is
defined. See~\cite{HindrySilvermanDGI,LangDG,SilvermanAEC} for this and
other basic properties of height functions.

\begin{remark}
\label{remark:minimalabsolutevalue}
Let~$\gp$ be a prime ideal of~$K$ corresponding to a place $v\in M_K$.
Then the largest value less than~$1$ of the normalized absolute value
attached to~$v$ is
\[
  \max_{\a\in\gp} |\a|_v  = q_\gp^{-1/[K_v:\QQ_v]}.
%%   = (\Norm_{K/\QQ}\gp)^{-1/[K_v:\QQ_v]}.
\]
\end{remark}

We set the following additional notation.

\begin{notation}
\item[$A/K$]
an abelian variety of dimension $g$ defined over $K$.
\item[$\End (A)$]
the endomorphism ring of $A$ over $\Kbar$.
\item[$\Lcal$]
an ample symmetric line bundle on $A/K$.
\item[$\hhat$]
the canonical height  $\hhat:A(\Kbar)\to[0,\infty)$ associated to~$\Lcal$.
If we need to make the dependence on~$A$ and/or~$\Lcal$ explicit, we will
write $\hhat_{\Lcal}$ or $\hhat_{A,\Lcal}$.
\item[$\Scal$]
a finite set of places of $K$ that includes at least all archimedean
places and all places of bad reduction of~$A$. From time to time we may
expand the set~$\Scal$, but it will depend only on~$A/K$, and not on
particular points of~$A$ or the extensions over which they are defined.
\end{notation}

%%%%%%%%%%%%%%%%%%%%%%%%%%%%%%%%%%%%%%%%%%%%%%%%%%%%%%%%%%%%%%%%%%%%%%
\section{A reduction step}
\label{section:reductionstep}

Theorem~\ref{theorem:maintheorem} deals with arbitrary abelian
varieties and ample line bundles. We begin with a reduction step
to a simpler situation.

\begin{proposition}[Reduction Step]
\label{proposition:reductionstep}
\hfil\break
\begin{parts}
\Part{(a)}
If Theorem~\ref{theorem:maintheorem} is true for geometrically simple
abelian varieties and very ample symmetric line bundles, then it is true for
all abelian varieties and all ample symmetric line bundles.
\Part{(b)}
If Theorem~\ref{theorem:maintheorem} is true for~$A/K'$ for some
finite extension~$K'/K$, then it is true for~$A/K$. In particular,
if~$A$ has complex multplication, it suffices to prove
Theorem~\ref{theorem:maintheorem} under the assumption 
that~$\End_K(A)=\End_\Kbar(A)$ and~$\End_K(A)\otimes\QQ\subset K$.
%%%%%%%%%
% Could also assume $\End(A)$ is integrally closed by taking an isogeny,
% but doesn't seem to be needed for the proof.
%%%%%%%%%
\end{parts}
\end{proposition}
\begin{proof}
(a)
Let $A/K$ be an arbitrary abelian variety and~$\Lcal$ an ample
symmetric line bundle on~$A/K$.
\par
Over an algebraically closed field, Poincar\'e's complete
reducibility theorem~\cite[Section~19]{Mumford} says that every
abelian variety decomposes, up to isogeny, into a product of simple
abelian varieties. We may thus find a finite extension~$K'$ and
geometrically simple abelian varieties~$A_1,\ldots,A_r$ defined
over~$K'$ so that there are isogenies
\[
  \f: A_1\times\cdots\times A_r\longrightarrow A
  \quad\text{and}\quad
  \psi: A \longrightarrow A_1\times\cdots\times A_r
\]
defined over~$K'$ with the property that $\f\circ\psi=[m]$ is
multiplication-by-$m$ on~$A$ for some integer~$m\ge1$.
\par
The line bundle~$\f^*\Lcal$ is ample on the 
product~$A_1\times\cdots\times A_r$, since~$\f$ is a finite 
morphism~\cite[Exercise~III.5.7(d)]{Hartshorne}.
We fix an integer~$n\ge1$ so that~$\f^*\Lcal^{\otimes n}$ is very
ample on the product, and we let~$\Lcal'_i=\f^*\Lcal^{\otimes n}\big|_{A_i}$
be the restriction to the~$i^{\text{th}}$~factor.
Then~$\Lcal'_i$ is a very ample line bundle on~$A_i$.
\par
Let $P\in A(\Kbar)$ and write $\psi(P)=(P_1,\ldots,P_r)$. Then
standard transformation properties of canonical heights (see
\cite[Chapters~4--5]{LangDG} or
\cite[Sections~B.3--B.5]{HindrySilvermanDGI}) allow us to compute
\begin{align*}
  \hhat_{A,\Lcal}(P)
 &= \frac{1}{m^2}\hhat_{A,\Lcal}([m]P) \\
 &= \frac{1}{m^2n}\hhat_{A,\Lcal^{\otimes n}}([m]P) \\
 &= \frac{1}{m^2n}
    \hhat_{A_1\times\cdots\times A_r,\f^*\Lcal^{\otimes n}}(\psi(P)) \\
 &= \frac{1}{m^2n}\left( \hhat_{A_1,\Lcal'_1}(P_1)
       + \hhat_{A_2,\Lcal'_2}(P_2) 
       + \cdots + \hhat_{A_r,\Lcal'_r}(P_r)\right).
\end{align*}
\par
By assumption, Theorem~\ref{theorem:maintheorem} is true
for each abelian variety~$A_i/K'$ and line bundle~$\Lcal'_i$, say
\[
  \hhat_{A_i,\Lcal'_i}(P)\ge C_i > 0
  \qquad\text{for all nontorsion $P\in A_i({K'}^\ab)$.}
\]
Since $K^\ab\subset {K'}^\ab$ and since
\[
  P\in A_\tors \Longleftrightarrow
  P_i\in (A_i)_\tors~\text{for all $1\le i\le r$},
\]
we conclude that if $P\in A(K^\ab)$ is a nontorsion point, then
\[
  \hhat_{A,\Lcal}(P) \ge \frac{1}{m^2n}\min_{1\le i\le r} C_i.
\]
This completes the proof of reduction step~(a).
\par
(b)
It is clear that Theorem~\ref{theorem:maintheorem} for~$K'$ implies
Theorem~\ref{theorem:maintheorem} for~$K$,
since~$K^\ab\subset{K'}^\ab$. The second part of~(b) follows, since
there is always a finite extension~$K'/K$ satisfying
$\End_{K'}(A)=\End_\Kbar(A)$ and
\text{$\End_{K'}(A)\otimes\QQ\subset {K'}$}.
\end{proof}

In view of Proposition~\ref{proposition:reductionstep}(a), we will
assume henceforth that the abelian variety~$A/K$ is geometrically simple
and that the line bundle~$\Lcal$ on~$A/K$ is very ample. We fix an embedding
\[
  \psi:A\longrightarrow \PP^n_K
\]
with $\psi^*(\OO(1))\cong\Lcal$. We fix a proper model~$\Acal/\O_K$
for~$A/K$ by letting $\Acal$ be the Zariski closure of $\psi(A)$ in
$\PP^n_{\O_K}$.  In particular, the embedding $\psi$ extends to a
closed immersion 
\[
  \Acal\longrightarrow \PP^n_{\O_K}.
\]
Adding a finite number of primes to the set~$\Scal$, we may assume
that~$\Acal$ is a group scheme over~$\Spec\Ocal_{K,\Scal}$ with
generic fiber~$A/K$. 
In particular, the group law on~$A$
commutes with reduction modulo~$\gp$  for $\gp\notin\Scal$. 
More precisely, if $P,Q,R\in A(K')$ for some finite extension $K'/K$,
and if~$\gP\in\Spec\O_{K'}$ is a prime of~$K'$ that does not lie
above a prime in~$\Scal$, then
\[
  P\equiv Q\pmod{\gP}
  \Longleftrightarrow
  P+R\equiv Q+R\pmod{\gP}.
\]

%%%%%%%%%%%%%%%%%%%%%%%%%%%%%%%%%%%%%%%%%%%%%%%%%%%%%%%%%%%%%%%%%%%%%%
\section{An elementary local-global height inequality}
\label{section:localglobalhtinequality}
For each place~$v \in M_K$, we define a logarithmic $v$-adic distance
function on~$\PP^n(K)$ by the formula
\[
  \d_v(x,y)
   = -\log\left(\frac{\displaystyle\Bigl(\max_{i,j} |x_iy_j-x_jy_i|_v\Bigr)}
          {\displaystyle\Bigl(\max_k|x_k|_v\Bigr)\Bigl(\max_k|y_k|_v\Bigr)}
          \right).
\]
Note that~$\d_v(x,y)$ does not depend on the choice of projective
coordinates for~$x$ and~$y$. The intuition is that
\[
  \d_v(x,y) = -\log|\text{$v$-adic distance from $x$ to $y$}|,
\]
so~$\d_v(x,y)$ becomes large as~$x$ and~$y$ become $v$-adically close
to one another. We set $\d_v(x,x)=\infty$ by convention.

\begin{lemma}
\label{lemma:heightvsdistance}
Let $K$ be a number field and let $T\subset M_K$ be any set of places
of~$K$. Then for all \emph{distinct} points $x,y\in\PP^n(K)$, 
\[
  h(x) + h(y) \geq \sum_{v \in T} n_v \d_v(x,y) - \log 2.
\]
\end{lemma}

\begin{proof}
For each $v\in M_K$, let $\e_v=1$ if~$v$ is archimedean, and $\e_v=0$
otherwise. The triangle inequality gives
\begin{align*}
  \max_{i,j} |x_i y_j - x_j y_i|_v 
  &\leq \max_{i,j} 2^{\epsilon_v} \max\bigl\{|x_iy_j|_v,|x_jy_i|_v\bigr\} \\
  &= 2^{\epsilon_v} \max_i \{ |x_i|_v \} \max_j \{ |y_j|_v \}.
\end{align*}
Hence
\begin{equation}
\label{equation:triangleinequality}
  \text{$\d_v(x,y)\ge -\e_v\log2$ for all $x,y\in\PP^n(K)$ and all
  $v\in M_K$.}
\end{equation}
\par
We are given that~$x\ne y$, so there exist indices~$\a,\b$ so that
\[
  x_\a y_\b \ne x_\b y_\a.
\]
We apply the product formula and the
inequality~\eqref{equation:triangleinequality} to compute
{\allowdisplaybreaks
\begin{align}
  h(x) + h(y) 
  &= \sum_{v\in M_K} 
      n_v\bigl(\log\max_k|x_k|_v+\log\max_k|y_k|_v\bigr) \notag\\
  &= \sum_{v\in M_K} 
      n_v\bigl(\max_{i,j}\log|x_iy_j-x_jy_i|_v + \d_v(x,y)\bigr) \notag\\
  &\geq \sum_{v\in M_K} 
      n_v\bigl(\log|x_\a y_\b-x_\b y_\a|_v + \d_v(x,y)\bigr) \notag\\
  &= \sum_{v\in M_K}  n_v\d_v(x,y)
     \qquad\text{(by the product formula)} \notag\\
  &\geq \sum_{v\in T} n_v\d_v(x,y) - \sum_{v\notin T} n_v\e_v\log2
      \qquad\text{(from \eqref{equation:triangleinequality})} \notag\\
  &\geq \sum_{v\in T} n_v\d_v(x,y) - (\log 2)\sum_{v\in M_K^\infty} n_v
      \qquad\text{(from def of $\e_v$)} \notag\\
  &\geq \sum_{v\in T}n_v\d_v(x,y) - \log 2.  \tag*{\qedsymbol}
\end{align}
}
\renewcommand{\qedsymbol}{}
\end{proof}

\begin{remark}
It is possible to put Lemma~\ref{lemma:heightvsdistance} into a more
general context, albeit with an inexplicit constant, using functorial
properties of Weil height functions associated to
subschemes. We briefly sketch the idea; see~\cite{SilvermanArithDistFunc} for the general
theory. Let~$H$ be a hyperplane in~$\PP^n$, and for
each $i\ne j$, let $D_{ij}\in\Div(\PP^n\times\PP^n)$ be the divisor of
type~$(1,1)$ defined by~$x_iy_j=x_jy_i$. Then
\[
  D_{ij}\sim\pi_1^*H+\pi_2^*H,\quad\text{so}\quad
  h_{D_{ij}} = h_H\circ\pi_1 + h_H\circ\pi_2 + O(1).
\]
The diagonal $\D\subset\PP^n\times\PP^n$ is the intersection $\D=\cap
D_{ij}$, hence~\cite[Theorem~4.1(b)]{SilvermanArithDistFunc} says that
\[
  h_\D \le \min_{i\ne j} h_{D_{ij}} + O(1) 
    = h_H\circ\pi_1 + h_H\circ\pi_2 + O(1).
\]
Finally we use the fact that~$h_\D$ is defined as a sum of local heights
\[
  h_\D = \sum_v \lambda_{\Delta,v},
\]
and the local heights are uniformly bounded below away from the 
support of~$\Delta$. Thus for any set of places~$T$,
\[
  h_\D \ge \sum_{v\in T} n_v\lambda_{\Delta,v} + O(1).
\]
Combining these inequalities and using the fact that
$\lambda_{\Delta,v}=\d_v+O(1)$ yields
\[
  h(x)+h(y) \ge \sum_{v\in T} n_v\d_v(x,y) + O(1)
  \quad\text{for all $(x,y)\in(\PP^n\times\PP^n)(\Kbar)\setminus\D$,}
\]
which is Lemma~\ref{lemma:heightvsdistance} with an inexplicit constant.
\end{remark}

We will apply Lemma~\ref{lemma:heightvsdistance} to our abelian
variety as follows.

\begin{proposition}
\label{proposition:primelowerbound}
There is a constant~$B=B(K,A,\Lcal)$ so that for all finite
extensions~$L/K$ and all pairs of distinct points~$P,Q\in A(L)$,
the following two inequalities are valid.
\begin{parts}
\Part{(a)}
Let $T\subset M_K$ be any set of places of~$K$. Then
\[
  \hhat(P) + \hhat(Q)
  \ge \sum_{v\in T} n_v\d_v(P,Q) - B.
\]
\Part{(b)}
Let~$\gP_1,\ldots,\gP_r$ be
primes of~$\O_{L}$ lying over primes of~$\O_K$ that are not 
in the set of excluded primes~$\Scal$. Assume further that
\[
  P\equiv Q\pmod{\gP_i}\qquad\text{for all $1\le i\le r$}.
\]
Then
\[
  \hhat(P) + \hhat(Q) 
  \ge \frac{1}{[L:\QQ]} \sum_{1\le i\le r} \log \Norm_{L/\QQ}\gP_i - B.
\]
\end{parts}
\end{proposition}
\begin{proof}
Recall that we have fixed an embedding $\psi:A\to\PP^n$ corresponding
to the line bundle~$\Lcal$. We apply
Lemma~\ref{lemma:heightvsdistance} to the points~$\psi(P)$
and~$\psi(Q)$ in~$\PP^n(L)$ and to the set of places~$v\in T$.  This
yields
\begin{equation}
  \label{equation:psiPQ}
  h(\psi(P))+h(\psi(Q)) 
   \ge \sum_{v\in T} n_v\d_v\bigl(\psi(P),\psi(Q)\bigr) - \log2.
\end{equation}
\par 
For all points $R\in A(\Kbar)$, functorial properties of height
functions give
\begin{align}
  h_{\PP^n,\OO(1)}(\psi(R))
  &= h_{A,\psi^*\OO(1)}(R) + O(1) \notag\\
  &= h_{A,\Lcal}(R) + O(1) 
  = \hhat_{A,\Lcal}(R) + O(1). 
  \label{equation:htonPn}
\end{align}
Combining~\eqref{equation:psiPQ} and~\eqref{equation:htonPn} gives
part~(a) of the proposition.
\par
Let $\{v_1,\ldots,v_r\}$ be the places of~$L$ corresponding to the set
of primes $\{\gP_1,\ldots,\gP_r\}$.  The congruence  $P\equiv
Q\pmod{\gP_i}$ implies that $\psi(P)\equiv \psi(Q)\pmod{\gP_i}$ (remember
that the map~$\psi$ has good reduction at all primes not in~$\Scal$).
Hence from Remark~\ref{remark:minimalabsolutevalue} we see that
\[
  \d_{v_i}\bigl(\psi(P),\psi(Q)\bigr)
   \ge \frac{1}{[L_{\gP_i}:\QQ_{\gP_i}]}\log \Norm_{L/\QQ}\gP_i.
\]
Substituting this into~(a) gives
\begin{align}
  \hhat(P)+\hhat(Q)
   &\ge \sum_{1\le i\le r} 
      \frac{n_{v_i}}{[L_{\gP_i}:\QQ_{\gP_i}]}
        \log \Norm_{L/\QQ}\gP_i - B \notag\\
   &= \frac{1}{[L:\QQ]} \sum_{1\le i\le r} \log \Norm_{L/\QQ}\gP_i - B.
   \tag*{\qedsymbol}
\end{align}
\renewcommand{\qedsymbol}{}
\end{proof}

%%%%%%%%%%%%%%%%%%%%%%%%%%%%%%%%%%%%%%%%%%%%%%%%%%%%%%%%%%%%%%%%%%%%%%
\section{A mod $\gp$ annihilating polynomial}
The proof of Theorem~\ref{theorem:maintheorem} is by induction on the
amount of ramification in a given finite abelian extension of $K$. The key to
handling the unramified case is to note that there is an element of
the group ring~$\ZZ[G_{\Kbar/K}]$ that simultaneously
annihilates~$A(K^\ab)$ modulo every prime of~$K^\ab$ lying
above~$\gp$, an idea already exploited in~\cite{SilvermanECAbExt} for
the case of elliptic curves. In this section we describe the
analogous fact for abelian varieties.

\begin{theorem}
\label{theorem:characpolyoffrob}
Let~$K/\QQ$ be a number field, let $\gp$ be a prime of~$K$, and let
$q=\Norm_{K/\QQ}(\gp)$.  Let~$A/K$ be an abelian variety of
dimension~$g$ with good reduction at~$\gp$, and let
\[
  \F_\gp(X) = \det\bigl(X-\Frob_\gp\bigm|T_\ell(A)\bigr) \in \ZZ[X]
\]
be the characteristic polynomial of Frobenius at~$\gp$. (Here $\ell$
is any prime number different from the residue characteristic of
$\gp$, and $T_\ell(A)$ is the $\ell$-adic Tate module of $A$.)
\par
Fix a prime $\bar\gp$ of~$\Kbar$ lying over~$\gp$ and
let~$\s\in(\bar\gp,\Kbar/K)\subset G_{\Kbar/K}$ be in the associated
Frobenius conjugacy class.
\begin{parts}
\Part{(a)}
Write~$\F_\gp(X)$ as
\[
  \F_\gp(X) = \sum_{i=0}^{2g} a_iX^i = \prod_{j=1}^{2g} (X-\a_j).
\]
Then the roots of~$\F_\gp$ satisfy $|\a_j|=q^{1/2}$ and the coefficients
of~$\F_\gp$ satisfy $|a_i|\le(4q)^g$.
\Part{(b)}
For all $P\in A(\Kbar)$,
\[
  \F_\gp(\s)P \equiv O \pmod{\bar\gp}.
\] 
(This congruence is taking place on a model $\Acal$ for~$A$ whose
special fiber $\tilde{A}$ is smooth.)
% for example on the N\'eron model.)
\Part{(c)}
If $P\in A(\Kbar)$ satisfies $\F_\gp(\s)P=O$, then~$P$ is a torsion point.
\end{parts}
\end{theorem}
\begin{proof}
(a) The equality~$|\a_j|=q^{1/2}$ is the Riemann hypothesis
for~$A(\FF_\gp)$, which was proven originally by Weil,
see~\cite[IV.21,~Theorem~4]{Mumford}. The coefficient~$a_j$ is the
$j^{\text{th}}$~symmetric polynomial of the roots of~$\F_\gp$, so 
\[
 |a_j| \le \binom{2g}{j}\max|\a_i|^j 
  \le 4^gq^{j/2} \le (4q)^g.
\]
\par
(b)
When reduced modulo~$\bar\gp$, the element~$\s\in G_{\Kbar/K}$ acts as
the $q$-power Frobenius map $f_q\in\End(\tilde A)$. Further,
the map~$\F_\gp(f_q)$ annihilates $T_\ell(\tilde A)$,
since~$\F_\gp$ is the characteristic polynomial of~$f_q$ acting
on~$T_\ell(\tilde A)$ and the Cayley-Hamilton theorem tells us
that a linear transformation satisfies its own characteristic equation.
\par
We have the general fact that if $B$ is an abelian variety over a
field $k$, then the map
\begin{equation}
  \label{equation:injectivityofEnd(B)}
  \End(B) \longrightarrow \End(T_\ell(B))
\end{equation}
is injective; see~\cite[IV.19,~Theorem~3]{Mumford} for the stronger
result that $\End(B)\otimes\ZZ_\ell\into\End(T_\ell(B))$.  We can
prove the injectivity of~\eqref{equation:injectivityofEnd(B)}
directly by noting that if~$\f\in\End(B)$ induces the zero map on~$T_\ell(B)$,
then $\f(B[\ell^n])=0$ for all $n\ge1$. Hence~$\f$ factors through
the isogeny~$[\ell^n]$, say $\f=\psi_n\circ[\ell^n]$. This implies that
$\deg(\f)=\deg(\psi_n)\deg([\ell])^n$. Since this holds for all $n\ge1$,
and since~$\deg(\f)$ and~$\deg(\psi_n)$ are integers and $\deg([\ell])>1$,
it follows that~$\deg(\f)=0$, and hence that~$\f=0$.
\par
Thus the fact that~$\F_\gp(f_q)$ annihilates~$T_\ell(\tilde{A})$ implies
that~$\F_\gp(f_q)=0$ as an element of~$\End(\tilde{A})$. In other
words,
\[
  \F_\gp(f_q)Q = O \qquad\text{for all $Q\in \tilde A(\bar\FF_\gp)$.}
\]
Finally, using the fact that the reduction map commutes with
the group law on~$A$, we see that for any~$P\in A(\Kbar)$, 
\[
  \F_\gp(\s)P \equiv \F_\gp(f_q)P \MOD{\bar\gp}
  \equiv O \MOD{\bar\gp}.
\]
\par
(c)
Let $P\in A(\Kbar)$ satisfy $\F_\gp(\s)P=O$. Fix a finite Galois
extension~$L/K$ with~$P\in E(L)$, say of degree $m=[L{:}K]$.
Then~$\s^m=1$ in~$G_{L/K}$, so in particular, $\s^mP=P$. Let
\[
  r = \Resultant(\F_\gp(X),X^m-1)\in\ZZ.
\]
The complex roots of~$X^m-1$ have absolute value~$1$ and from~(a), the
complex roots of~$\F_\gp(X)$ have absolute value~$q^{1/2}$, so the two polynomials
have no complex roots in common. It follows that~$r\ne0$.
\par
The resultant of two polynomals in~$\ZZ[X]$ is an element of the ideal
that they generate~\cite[Chapter IV, Section 8]{LangAlg},
so we can find polynomials $u(X),v(X)\in\ZZ[X]$ satisfying
\[
  u(X)\F_\gp(X)+v(X)(X^m-1)=r.
\]
Substituting $X=\s$ gives the identity 
\[
  u(\s)\F_\gp(\s) + v(\s)(\s^m-1)=r
\]
in the group ring $\ZZ[G_{\Kbar/K}]$. Hence
\[
  rP = u(\s)\bigl(\F_\gp(\s)P\bigr)
      + v(\s)\bigl((\s^m-1)P\bigr) = O,
\]
so~$P$ is a point of finite order.
\end{proof}

%%%%%%%%%%%%%%%%%%%%%%%%%%%%%%%%%%%%%%%%%%%%%%%%%%%%%%%%%%%%%%%%%%%%%%
\section{Selection of a ``good'' prime}
In order to prove our main result (Theorem~\ref{theorem:maintheorem}),
we will show that there is a prime $\gp$ of $K$ and 
positive constants~$C_1,C_2,C_3$, depending
only on~$A/K$ and~$\Lcal$, such that
\[
  \hhat_{A,\Lcal}(P) 
  \geq C_1 \frac{\log\Norm_{K/\QQ}\gp - C_2}{(\Norm_{K/\QQ}\gp)^{C_3}}
  > 0
  \qquad\text{for all nontorsion $P\in A(K^\ab)$.}
\]
% for some prime~$\gp$ of~$K$ that makes the righthand side positive.
The next proposition will help us choose such a prime $\gp$.
We assume from now on that the abelian variety $A$ is $\Kbar$-simple and
that $\Lcal$ is very ample.
% and that $\End_K(A) = \End_{\Kbar}(A)$.
By Proposition~\ref{proposition:reductionstep}, these assumptions are harmless.

% \par
% For the statement of the next proposition, recall that an abelian
% variety $A$ of dimension $g$ over a field $k$ of characteristic $p>0$ is called {\em
% ordinary} if $\# A[p](\overline{k}) = p^g$.

\begin{proposition}
\label{proposition:primeconditions}
There exist infinitely many primes~$\gp$ of~$K$ satisfying 
the following conditions:
\begin{enumerate}
\item
$\gp$ is unramified of degree one, so that in particular $N^K_\QQ \gp$
is a prime number $p$. 
% Let \text{p=$\Norm_{K/\QQ}\gp$} denote the characteristic of
% the residue field~$\O_K/\gp$.
\item
Reduction modulo~$\gp$ commutes with the embedding $\psi:A\into\PP^n$
coming from the very ample symmetric line buncle~$\Lcal$.
\item
The abelian variety~$A$ has good reduction at~$\gp$, so that in particular,
reduction modulo~$\gp$ commutes with the group law on~$A$.  
\item
$p > \exp([K:\QQ](B+1))$, where~$B$ is the
constant appearing in Proposition~\ref{proposition:primelowerbound}.
\item
One of the following holds:
% The final condition on~$\gp$ splits into two cases:
\begin{enumerate}
\item
Either~$A(K^\ab)$ has no points of order~$p$, or
\item
$A$ has complex multiplication and ordinary reduction at~$\gp$.
\end{enumerate}
\end{enumerate}
\end{proposition}
\begin{proof}
It is clear that conditions~(2),~(3), and~(4) exclude only finitely
many primes. If~$A$ does not have complex multiplication, then the
same is true for condition~(5a) by a theorem of Zarhin
(Theorem~\ref{theorem:noptorsion} below).  On the other hand, if~$A$
does have complex multiplication, then it is well-known that~$A$ has
ordinary reduction at all sufficiently large primes which split
completely in some fixed finite extension of $K$ (see Theorem~\ref{theorem:ordinaryCMtheorem} below).
% Hence for CM abelian varieties, condition~(5b) holds for all
% primes that split completely in some fixed finite extension of~$K$.
Since there are infinitely many rational primes which split completely
in any given number field, this shows that we can find (infinitely many) primes
satisfying all five conditions.
\end{proof}

\begin{theorem}
\label{theorem:ordinaryCMtheorem}
Let $A$ be a $g$-dimensional abelian 
variety over the number field $K$, 
and suppose that $L := (\End  A)\otimes \QQ$ 
is a field of degree $2g$ over $\QQ$.  Let   
$\gp$ be a degree-$1$ prime of $K$ and let $p$ be the rational
prime under $\gp$.  Suppose that
\begin{enumerate}
\item $A$ has good reduction at $\gp$; 
\item Every endomorphism of $A$ has good reduction at $\gp$;
\item $\End_K(A) = \End_{\Kbar}(A)$; and 
\item $p$ is unramified in $L$. 
\end{enumerate}
Then $A$ has ordinary reduction at $\gp$. 
\end{theorem}

\begin{proof}
Since the property of being ordinary is preserved by isogenies and
products, we may assume without loss of generality that $A$ is simple
and that $A$ is {\em principal}, i.e., that
$\End(A)$ is the maximal order in $(\End A)\otimes \QQ$.
It then follows from \cite[Chapter~III, Theorem~2]{Shimura} and our hypotheses
that the reduction~$\tilde A$ of~$A$ modulo~$\gp$ is also simple, and that the natural
reduction map gives an \emph{isomorphism}
\[
  \End(A)\otimes\QQ 
  \stackrel{\sim}{\longrightarrow} 
  \End(\tilde A)\otimes\QQ.
\]
In particular, since $L := \End(A)\otimes\QQ$ is a CM field of degree $2g$, we conclude that
$\End(\tilde A)\otimes\QQ = L$ is a CM field of degree $2g$ as well.

Let $\pi \in L$ be the Frobenius morphism of $\tilde A$ over $\FF_p$; 
then $\pi$ is in fact an element of the ring of integers $\O_L$ of $L$. Let
$\overline{\pi}$ be the complex conjugate of $\pi$.  A consideration of degrees and
the fact that $L$ is a CM field shows that $\pi \overline{\pi} = p$.  Since $[L:\QQ] =
2g$, it follows from \cite[Theorem 2, page 140]{Tate}
that $L = \QQ(\pi)$ and that $\F$, the characteristic polynomial of
$\pi$ acting on $T_\ell(\tilde A)$ for some prime $\ell \neq p$, has no
multiple roots.  Since every conjugate of $\pi$ satisfies the polynomial
$\F$ as an endomorphism of $\tilde A$, it follows that 
\[ 
\F(X) = \prod_{i=1}^{2g} (X - \pi_i), 
\] 
where $\pi_1,\ldots,\pi_{2g}$ are the conjugates of $\pi$.
 
Since $p$ is unramified in $\O_L$, it follows from the relation $p = \pi
\pi'$ that for any prime ideal $\gq$ of $\O_L$ lying over $p$ and
for each $i$, exactly one of $\pi_i$ and $\overline{\pi}_i$ is divisible by $\gq$.
% exactly $g$ of the $2g$ roots of $\F$ are divisible
% by any fixed prime ideal of $R$ lying over $p$.  
This implies that $\tilde A$ is ordinary (see \cite[proof of Proposition 7.1]{Waterhouse}).
\end{proof}

\begin{theorem}
\label{theorem:noptorsion}
Let $A/K$ be a geometrically simple abelian variety that does not have complex
multiplication (over~$\Kbar$). Then~$A(K^\ab)_\tors$ is finite.
In particular,
\[
  A(K^\ab)[p] = 0\quad\text{for all but finitely many primes~$p$.}
\]
\end{theorem}
\begin{proof}
This is proven by Zarhin in~\cite{Zarhin}, using methods developed by
Faltings~\cite{Faltings} in his proof of Tate's isogeny conjecture.
See also~\cite{Ruppert,Serre,SerreImageOfGalois,SerreImageOfGaloisAV}.
\end{proof}

\begin{remark}
\label{remark:Ruppert}
In the statament of Theorem~\ref{theorem:noptorsion}, one can replace
the hypothesis that $A/K$ is a geometrically simple abelian variety
without complex multiplication over~$\Kbar$ by the hypothesis that
$A/K$ is an abelian variety having no abelian subvariety with
complex multiplication over $K$ (see \cite{Ruppert}).
\end{remark}

For the remainder of this paper, we fix a prime~$\gp$ of~$K$
satisfying the conditions described in
Proposition~\ref{proposition:primeconditions}, and we let
$p=\Norm_{K/\QQ}\gp$.

%%%%%%%%%%%%%%%%%%%%%%%%%%%%%%%%%%%%%%%%%%%%%%%%%%%%%%%%%%%%%%%%%%%%%%
\section{Proof of the main theorem in the unramified case}
Let~$A/K$ be an abelian variety and~$\Lcal$ a line bundle on~$A/K$.
As we have already mentioned, we may assume 
% In order to prove Theorem~\ref{theorem:maintheorem}, it suffices by
% Proposition~\ref{proposition:reductionstep} to consider the case
that~$A$ is geometrically simple and~$\Lcal$ is very ample.  Recall
that we have fixed a prime~$\gp$ of~$K$ satisfying the conditions in
Proposition~\ref{proposition:primeconditions}, and that $p$ is the
residue characteristic of $\gp$.
\par
Let $P\in A(K^\ab)$ be a nontorsion point.
The proof of Theorem~\ref{theorem:maintheorem} is by induction
on $\ord_\gp(\gf_\gp(K(P)/K))$, where $\gf_\gp(L/K)$ denotes the
local conductor of the abelian extension~$L/K$
at~$\gp$.  For our purposes, we define $\gf_\gp(L/K)$ to be the
smallest positive integer $m$ such that $L_\gq \subseteq
\QQ_p(\zeta_m)$, where $\gq$ is a prime of $L$ lying over $\gp$.  It
follows from local class field theory that $\gf_\gp(L/K)$ exists and
is well-defined, and that $\ord_\gp(\gf_\gp(L/K)) \geq 1$ if and only
if $L/K$ is ramified at $\gp$.
% The proof of Theorem~\ref{theorem:maintheorem} is by induction
% on the ramification index of~$K(P)/K$ at~$\gp$. 
In this section we begin the
induction by proving the unramified case. 

\begin{theorem}
\label{theorem:maintheoremunramifiedcase}
Let $L\subset K^\ab$ be unramified at~$\gp$ and let~$P\in A(L)$ be
a nontorsion point. Then
\[
  \hhat(P) \ge \frac{1}{(12p)^{2g}}.
\]
(Note that this gives a lower bound for~$\hhat(P)$ that is independent
of~$L$ and~$P$.)
% since $p$ and $\Norm_{K/\QQ}\gp$ are fixed via
% Proposition~\ref{proposition:primeconditions} before~$L$ and~$P$ are
% selected.)
\end{theorem}
\begin{proof}
Factor the prime~$\gp$ in~$L$ as
\[
  \gp\O_L = \gP_1\gP_2\cdots\gP_r.
\]
Since $\gp$ is unramified in $L$, the $\gp_i$'s are all distinct.  Moreover,
the fact that the extension~$L/K$ is abelian implies that all of the
primes $\gP_1,\ldots,\gP_r$ have the same associated Frobenius
element in~$G_{L/K}$, say
\[
  \s = (\gP_i,L/K) \in G_{L/K}\qquad\text{for all $1\le i\le r$.}
\]
\par
Let $\F_\gp(X)\in\ZZ[X]$ be the characteristic polynomial of Frobenius
acting on~$T_\ell(A)$. Theorem~\ref{theorem:characpolyoffrob}(b)
tells us that
\[
  \F_\gp(\s)P\equiv O\MOD{\gP_i}
  \qquad\text{for all $1\le i\le r$}.
\]
% Note how we are using here the fact that every~$\gP_i$ has
% the same Frobenius element in~$G_{L/K}$, which in turn is a consequence
% of the fact that~$L/K$ is abelian.
\par
Theorem~\ref{theorem:characpolyoffrob}(c) and our assumption
that~$P$ is a nontorsion point tell us that~$\F_\gp(\s)P\ne O$. 
Hence we can apply Proposition~\ref{proposition:primelowerbound}(b)
to the \emph{distinct} points~$\F_\gp(\s)P$ and~$O$. Since~$\hhat(O)=0$,
this yields
\[
  \hhat(\F_\gp(\s)P) 
  \ge \frac{1}{[L:\QQ]}\sum_{i=1}^r \log\Norm_{L/\QQ}\gP_i - B
\]
for a constant~$B$ that is independent of~$L$ and~$P$.
The factorization $\gp=\gP_1\cdots\gP_r$ implies that
\[
  \sum_{i=1}^r \log\Norm_{L/\QQ}\gP_i 
  = \log\Norm_{L/\QQ}\gp
  = [L:K]\log\Norm_{K/\QQ}\gp,
\]
so we obtain the lower bound
\begin{equation}
  \label{equation:htfroblowerbd}
  \hhat(\F_\gp(\s)P)
  \ge\frac{\log p}{[K:\QQ]} - B.
\end{equation}
(Remember that~$\gp$ has degree one over~$\QQ$.)
\par
Next we write $\F_\gp(X)=\sum a_jX^j$ and compute
\begin{align*}
  \hhat(&\F_\gp(\s)P)
  =\hhat\biggl( \sum_{j=0}^{2g} a_j\s^jP \biggr) \\
  &\le (2g+1) \sum_{j=0}^{2g} \hhat(a_j\s^jP) 
    \qquad\text{parallelogram law (see Remark~\ref{remark:parallelogramlaw}),} \\
  &= (2g+1) \sum_{j=0}^{2g} a_j^2\hhat(\s^jP)
     \qquad\text{since $\hhat$ is a quadratic form,} \\
  &\le  (2g+1) \sum_{j=0}^{2g} (4p)^{2g} \hhat(\s^jP)
      \qquad\text{from Theorem~\ref{theorem:characpolyoffrob}(a),} \\
  &= (2g+1)^2(4p)^{2g}\hhat(P)
      \qquad\text{since $\hhat$ is Galois invariant.}
\end{align*}
Combining this with~\eqref{equation:htfroblowerbd} and the trivial
estimate \text{$2g+1\le 3^g$} gives
\[
  \hhat(P) \ge \frac{1}{(12p)^{2g}}
     \left(\frac{\log p}{[K:\QQ]} - B\right).
\]
Finally, we recall that condition~(4) of
Proposition~\ref{proposition:primeconditions} says that
% $p$ satisfies
$p>\exp\bigl([K:\QQ](B+1)\bigr)$, which yields the stated lower bound
$\hhat(P) \ge 1/(12p)^{2g}$.
\end{proof}

\begin{remark}
\label{remark:parallelogramlaw}
During the proof of Theorem~\ref{theorem:maintheoremunramifiedcase} we
made use of the following generalized parallelogram law.  For any
quadratic form~$Q$, it is easy to check the formal identity
\[
  Q\biggl(\sum_{i=1}^t x_i\biggr) + \frac{1}{2}\sum_{i,j=1}^t Q(x_i-x_j)
  = t\sum_{i=1}^t Q(x_i).
\]
Hence if~$Q$ is positive semidefinite, then
\[
  Q\biggl(\sum_{i=1}^t x_i\biggr) \le t\sum_{i=1}^t Q(x_i).
\]
\end{remark}

\begin{remark}
The proof of Theorem~\ref{theorem:maintheoremunramifiedcase} actually
shows that
\[
  \hhat(P) \ge \frac{1}{(2g+1)^2(4p)^{2g}}
     \left(\frac{\log p}{e_\gp(L/K)[K:\QQ]} - B\right),
\]
where~$e_\gp(L/K)$ is the ramification index of~$L/K$. This suffices
to prove our main theorem for extensions whose ramification index
at~$\gp$ is bounded, but it will not handle highly ramified
extensions.
\end{remark}

%%%%%%%%%%%%%%%%%%%%%%%%%%%%%%%%%%%%%%%%%%%%%%%%%%%%%%%%%%%%%%%%%%%%%%
\section{A generalized Amoroso-Dvornicich lemma}
%\section{Preliminaries for ramified abelian extensions}

In order to deal with the ramified case of Theorem~\ref{theorem:maintheorem},
we recall the following lemma of Amoroso and Dvornicich
(see~\cite[Lemma~2]{AmorosoDvornicich}, and also
also~\cite{AmorosoZannier,Baker,SilvermanECAbExt}). We refer the reader
to~\cite{SilvermanECAbExt} for a proof of this particular formulation.

\begin{lemma}[Amoroso-Dvornicich \cite{AmorosoDvornicich}]
\label{lemma:inertiacongruence}
Let $K/\QQ$ be a number field, 
let~$\gp$ be a degree~1 prime of~$K$ with residue characteristic~$p$,
% Maybe say: a prime in $K$ such that $K_\p \cong \QQ_p$??
and let~$L/K$ be an abelian extension
that is ramified at~$\gp$. Let~$\gP$ be a prime of~$L$ lying
over~$\gp$, let~$\Ocal_{L,\gP}$ denote the localization of~$L$
at~$\gP$, and let $I_{L/K}$ be the inertia group at $\gP$.  
Then there exists an element~$\t\in I_{L/K}$ with~$\t\ne1$
such that
\[
  \t(\a)^p \equiv \a^p \pmod{p\Ocal_{L,\gP}}
  \qquad\text{for all $\a\in\Ocal_{L,\gP}$.}
\]
(Note that the strength of this result is that the congruence is
modulo~$p$, and not merely modulo~$\gP$.)
\end{lemma}

\begin{remark}
\label{remark:inertiacongruenceremark}
Suppose that $L_\gP=\QQ_p(\z_m)$ for some integer~$m$ divisible by~$p$. Then the
proof of Lemma~\ref{lemma:inertiacongruence} shows that we may take~$\t$ to be
any nontrivial element of $\Gal\bigl(\QQ_p(\z_m)/\QQ_p(\z_{m/p})\bigr)$,
considered as a subgroup of
\[
  \Gal(\QQ_p(\z_m)/\QQ_p)  
  =\Gal(L_\gP/K_\gp) 
  \subset\Gal(L/K).
\]
\end{remark}

We now prove a version of the lemma of Amoroso and Dvornicich that
applies to varieties and maps with a particular type of inseparable reduction.

\begin{proposition}
\label{proposition:ramifieddistancebound}
Let $K/\QQ$ be a number field and let~$\gp$ be a degree~1 prime
of~$K$.  Let $X/K\subset\PP^n_K$ be a 
% nonsingular 
variety and let $\f:X\to X$ be a finite $K$-morphism.  Fix a model
$\Xcal/{\Ocal_{K,\gp}}\subset\PP^n_{\Ocal_{K,\gp}}$ and
let $\F:\Xcal\to\Xcal$ denote the extension of~$\f$
to~$\Xcal$. Make the following assumptions:
\begin{enumerate}
\item
The scheme $\Xcal/{\Ocal_{K,\gp}}$ is embedded as a 
% smooth and 
proper ${\Ocal_{K,\gp}}$-subscheme of~$\PP^n_{\Ocal_{K,\gp}}$.
\item
The map $\F:\Xcal\to\Xcal$ is a finite ${\Ocal_{K,\gp}}$-morphism.
\item
The restriction of~$\F$ to the special fiber, i.e., the reduction 
$\widetilde{\F} : \widetilde{\Xcal} \to \widetilde{\Xcal}$,
of~$\F$ modulo~$\gp$, 
factors through the Frobenius map ${\rm Frob} : \widetilde{\Xcal} \to \widetilde{\Xcal}^{(p)}$.
% is \emph{inseparable}.
\end{enumerate}
\par
Let $L/K$ be an abelian extension which is ramified at $\gp$, let $\gP|\gp$, and 
let $\t\in I_{L/K}$ be as in Lemma~\ref{lemma:inertiacongruence}.
Then for all points $P\in X(L)$,
\begin{equation}
  \label{equation:ramifieddistancebound}
  \d_\gP\bigl(\f(\t(P)),\f(P)\bigr) \ge \log p.
\end{equation}
(We recall that~$\d_\gP$ is the (logarithmic) $\gP$-adic distance
function defined in section~\ref{section:localglobalhtinequality}.)
\end{proposition}

% This remark is wrong and should be deleted.
%
% \begin{remark}
% \label{remark:inseparable}
% Recall that a dominant rational map $\psi : X \to Y$ of varieties over
% a field $k$ is called {\em inseparable} if the induced map $\psi^* :
% k(Y) \to k(X)$ on function fields is not separable.  By field theory,
% if $k$ is perfect, $\psi$ is inseparable, and $f \in k(X)$, then $\psi^*f = f\circ
% \psi$ is a $p$th power in $k(Y)$.  Also, by
% \cite[Chapter III, Theorem 1]{Shafarevich}, if $\psi^* :
% H^0(Y,\Omega^1) \to H^0(X,\Omega^1)$ is not injective then $\psi$ is inseparable.
% \end{remark}

\begin{remark}
The crucial property of the
estimate~\eqref{equation:ramifieddistancebound} is that the lower
bound for the $\gP$-adic distance does not depend on the ramification
degree of~$\gP$.  Since~$\t$ is in the inertia group, we know that
$\t(P)\equiv P\MOD{\gP}$. Hence
Remark~\ref{remark:minimalabsolutevalue} gives the trivial lower bound
\[
  \d_\gP(\t(P),P) 
  \ge \frac{1}{[L_\gP:\QQ_p]} \log\Norm_{L/\QQ}\gP 
  = \frac{\log p}{e_{L/\QQ}(\gP)}.
\]
This would not suffice for our purposes.
\end{remark}

\begin{proof}
Without loss of generality, we may replace~$K$ and~$L$ by~$K_\gp$
and~$L_\gP$, respectively. 
\par
Let~$\Ucal/\Ocal_K$ and~$\Vcal/\Ocal_K$ be affine open subsets
of~$\Xcal/\Ocal_K$ with~$\F(\Ucal)\subset\Vcal$. Choose
affine coordinates (i.e., generators of the affine coordinate ring
as an $\Ocal_K$-algebra) $\bfx=(x_1,\ldots,x_r)$ on~$\Ucal$ and
simlarly~$\bfy=(y_1,\ldots,y_s)$ on~$\Vcal$.
The map~$\F:\Ucal\to\Vcal$ is given by polynomials
\[
  \F^*\bfy = \bfA(\bfx) = (A_1(\bfx),\ldots,A_s(\bfx))
  \qquad\text{with $A_1,\ldots,A_s \in\Ocal_K[\bfx]$.}
\]
We are given that~$\F\bmod\gp$ factors through the Frobenius map. 
% Previous condition: the reduction of~$\F$ modulo~$\gp$ is inseparable.
It follows that
\begin{equation}
  \label{equation:formofinseparablemap}
  \F^*\bfy = \bfB(\bfx^p)+p\bfC(\bfx)
  \qquad\text{with $\bfB,\bfC\in\Ocal_K[\bfx]^s$,}
\end{equation}
where for notational convenience we write $\bfx^p=(x_1^p,\ldots,x_r^p)$.
(Remember that~$\gp$ has degree~1, so~$p$ is a uniformizer for~$\gp$.)
\par
Let $P\in X(L)=\Xcal(\Ocal_L)$.  Choose an open affine
subsets and local coordinates as above with~$P\in\Ucal$ and
$\F(P)\in\Vcal$. (More formally,~$P$ is really a morphism
$P:\Spec(\Ocal_L)\to\Xcal$ over~$\Spec(\Ocal_K)$, and
we choose~$\Ucal$ so that $P(\Spec\Ocal_L)\subset\Ucal$. Similarly
for~$\F(P)$ and~$\Vcal$.) We compute
\begin{align}
   \bfy(\F&(\t P))-\bfy(\F(P)) \notag\\
   &= \bigl(\bfB(\bfx(\t P)^p)+p\bfC(\bfx(\t P))\bigr)
     - \bigl(\bfB(\bfx(P)^p)+p\bfC(\bfx(P))\bigr) \notag\\
   &\omit\hfill\text{from \eqref{equation:formofinseparablemap},} \notag\\
   &= \bigl(\bfB(\bfx(P)^p + p\bfa)+p\bfC(\bfx(\t P))\bigr)
     - \bigl(\bfB(\bfx(P)^p)+p\bfC(\bfx(P))\bigr) \notag\\
    &\omit\hfill\text{from Lemma~\ref{lemma:inertiacongruence},
        where $\bfa\in\Ocal_L^r$}, \notag\\
   &\equiv 0 \pmod{p\Ocal_L}.
   \label{equation:F(tP)-F(P)}
\end{align}
\par
The $\gP$-adic distance between any two points~$Q_1,Q_2\in\Vcal(\Ocal_L)$
is given by
\begin{equation}
  \label{equation:affinedistance}
  \d_\gP(Q_1,Q_2) 
    = -\log \max_{1\le i\le s} \bigl|y_i(Q_1)-y_i(Q_2)\bigr|_\gP.
\end{equation}
(Note that~$\bfy$ gives affine coordinates on~$\Vcal$, so projective
coordinates are~$[1,y_1,\ldots,y_s]$.) We apply this formula
with~$Q_1=\F(\t P)$ and~$Q_2=\F(P)$ to obtain the desired result:
\begin{align*}
  \d_\gP(\F(\t P),\F(P))
  &= -\log \max_{1\le i\le s} \bigl|y_i(\F(\t P))-y_i(\F(P))\bigr|_\gP  
    \quad\text{from \eqref{equation:affinedistance},} \\
  &\ge -\log|p|_\gp
    \quad\text{from \eqref{equation:F(tP)-F(P)},} \\
  &=\log p.  \tag*{\qedsymbol}
\end{align*}
\renewcommand{\qedsymbol}{}
\end{proof}

In order to apply Proposition~\ref{proposition:ramifieddistancebound},
we will use the following result.

\begin{lemma}
\label{lemma:fontainelemma}
If $\widetilde{A}$ is an abelian variety over a perfect field $k$ of
characteristic $p>0$, then $[p] : \widetilde{A} \to \widetilde{A}$ factors
through the Frobenius map ${\rm Frob} : \widetilde{A} \to \widetilde{A}^{(p)}$.
\end{lemma}

\begin{proof}
By the theory of isogenies and quotients, it suffices to prove that
${\rm ker}({\rm Frob}) \subseteq {\rm ker}([p])$ as finite flat group schemes
over $k$, i.e., that the connected group scheme ${\rm ker}({\rm Frob})$ is
killed by $p$.  This holds for any formal commutative
group scheme $G$ over $k$, as follows from the theory of the
Verscheibung operator (see \cite[Chapter I, Section 7.5]{Fontaine} 
or \cite[Expos{\'e} $\rm{VII_A}$, Sections 4.2--4.3]{SGA3} for details).
% for example from the identity $[p]_G = {\rm Ver} \circ {\rm Frob}$ 
\end{proof}

%%%%%%%%%%%%%%%%%%%%%%%%%%%%%%%%%%%%%%%%%%%%%%%%%%%%%%%%%%%%%%%%%%%%%%

\section{Proof of the main theorem in the ramified, non-CM case}
We now consider the case in which the point $P\in A(K^\ab)$ is 
defined over a field that is ramified at~$\gp$.  In this section, we
assume furthermore that $A$ does not have complex multiplication.  We will use the
following induction hypothesis.
\begin{center}
   \Hypothesis($e$):\quad
   $\left(\parbox{.75\hsize}{\noindent
     For all fields $K'\subset K^\ab$ whose 
     $\gp$-rami\-fi\-ca\-tion index satisfies
     \text{$e_\gp(K'/K)\le e$} and all nontorsion points \text{$P\in A(K')$},
     the height of~$P$ satisfies \text{$\hhat(P) \ge \dfrac{1}{(12p)^{2g}}$}.
   }\right)$
\end{center}

We note that Theorem~\ref{theorem:maintheoremunramifiedcase} shows
that \Hypothesis(1) is true.

\begin{theorem}
\label{theorem:maintheoremramifiednonCMcase}
Suppose that $A$ does not have complex multiplication (over $\Kbar$).
% Let $L\subset K^\ab$ be ramified at~$\gp$ 
Let $L$ be a finite abelian extension of $K$
and assume that \Hypothesis($e$)
is true for all $e<e_\gp(L/K)$. Then 
% is true for all $f<\ord_\gp(\gf_\gp(L/K))$. Then 
\[
  \hhat(P) \ge \frac{1}{(12p)^{2g}} 
    \qquad\text{for all nontorsion $P\in A(L)$.}
\]
Hence by induction, \Hypothesis($e$) is true for all~$e\ge1$.
\end{theorem}

\begin{proof}
As noted above, we have already shown that~\Hypothesis(1) is true, so
we may assume the~$L/K$ is ramified at~$\gp$. Let $P\in A(L)$ be a
nontorsion point. If~$K(P)/K$ is less ramifed at~$\gp$ than~$L/K$
(i.e., if $e_\gp(K(P)/K) < e_\gp(L/K)$),
then we are done by the induction hypothesis, so we may assume that
\[
  e_\gp(K(P)/K) = e_\gp(L/K).
\]
Let $\t\in G_{L/K}$ be chosen as in
Lemma~\ref{lemma:inertiacongruence}.
\par
% The multiplication-by-$p$ map $[p]:A\to A$ acts on the space of
% invariant one-forms via $[p]^*\omega=p\omega$, hence in characteristic~$p$
% it annihilates~$H^0(A,\Omega^1)$. Thus~$[p]$ has inseparable reduction
% by Remark~\ref{remark:inseparable},
Using Lemma~\ref{lemma:fontainelemma}, we now apply Proposition~\ref{proposition:ramifieddistancebound}
to the point~$P\in A(L)$ and the map~$\F=[p]$. The proposition tells
us that
\[
  \d_\gP([p](\t P),[p](P))\ge\log p
  \qquad\text{for all primes~$\gP|\gp$.}
\]
\par
Summing over all of the primes dividing~$\gp$, we obtain
a lower bound that is independent of the ramification degree
of~$L/K$:
\begin{align}
  \sum_{\gP|\gp} \frac{[L_\gP:\QQ_p]}{[L:\QQ]} \d_\gP(([p](\t P),[p](P))
  &\ge
  \frac{[K_\gp:\QQ_p]}{[K:\QQ]}
  \sum_{\gP|\gp} \frac{[L_\gP:K_\gp]}{[L:K]} \log p \notag\\
  &=\frac{\log p}{[K:\QQ]}.
  \label{equation:ramifdeltalowerbd}
\end{align}
(Remember that~$\gp$ is a degree one prime, so $K_\gp=\QQ_p$.)
\par
Assume for the moment that $[p](\t P)\ne[p](P)$. Then we can apply
Proposition~\ref{proposition:primelowerbound}(a)
to the \emph{distinct} points~$[p](\t P)$ and~$[p](P)$, which
gives the lower bound
\begin{align*}
  \hhat([p](\t P)) + \hhat([p]P)
  &\ge \sum_{\gP|\gp} \frac{[L_\gP:\QQ_p]}{[L:\QQ]}
               \d_\gP(([p](\t P),[p](P)) - B \\
  &\ge \frac{\log p}{[K:\QQ]} - B 
    \qquad\text{from \eqref{equation:ramifdeltalowerbd},} \\
  &\ge 1
    \qquad\text{from Proposition~\ref{proposition:primeconditions}(4).} 
\end{align*}
\par
Using the quadratic nature and Galois invariance of the canonical
height yields
\[
  \hhat([p](\t P)) + \hhat([p]P) = 2p^2\hhat(P).
\]
Hence $\hhat(P) \ge 1/2p^2$, which is considerably stronger than the
desired result. This completes the proof under the assumption that
$[p](\t P)\ne[p](P)$.
\par
Finally, let $Q=\t P - P$, and suppose that~$[p](Q)=O$.  In
particular,~$Q\in A(L)_\tors$. 
% We now consider two cases, depending
% on whether our chosen prime~$\gp$ satisfies Condition~5(a) or
% Condition~5(b) of Proposition~\ref{proposition:primeconditions}.
% We begin with the easier case of~5(a),which says that~$A(K^\ab)$ 
% contains no nontrivial points of order~$p$. 
Since we are assuming that $A$ does not have complex multiplication, 
we may assume that our chosen prime~$\gp$ satisfies Condition~5(a)
of Proposition~\ref{proposition:primeconditions}, i.e., 
that~$A(K^\ab)$ contains no nontrivial points of order~$p$. 
This implies immediately
that~$Q=O$, and hence that~\text{$\t P=P$}. Therefore~$P$ is defined over
the fixed field~$L^\t$ of~$\t$. However,~$\t$ is a nontrivial element
of the inertia group at~$\gp$, so
\[
  e_\gp(L^\t/K) < e_\gp(L/K)\qquad\text{(strict inequality)}.
\]
Therefore~$\hhat(P)>1/(12p)^{2g}$ by the induction hypothesis applied to
the field~$L^\t$.
% \par
% It remains to deal with the case that~$A(K^\ab)$ may contain a large
% amount of torsion, which by Theorem~\ref{theorem:noptorsion}
% is essentially only when~$A$ has complex multiplication. 
% The moderately complicated argument for this case,
% and thus the completion of the proof of our
% main theorem, is given in the next section.
% \renewcommand{\qedsymbol}{}
\end{proof}

\begin{remark}
We could just as well have proved Theorem
by induction on the local conductor $\gf_\gp(L/K)$ (rather than on the
ramification index $e_\gp(L/K)$).  In the next section, we will
necessarily have to use the local conductor, in order to assume that
the field $L$ contains enough roots of unity.
\end{remark}

%%%%%%%%%%%%%%%%%%%%%%%%%%%%%%%%%%%%%%%%%%%%%%%%%%%%%%%%%%%%%%%%%%%%%%
\section{Completion of the proof --- $L/K$ is ramified and 
$A$ has complex multiplication}

An examination of the proof of Theorem~\ref{theorem:maintheoremramifiedCMcase}
given in the previous section 
shows that we have reduced the proof of Theorem~\ref{theorem:maintheorem} to the
case that~$A$ has complex multiplication, $A$~has ordinary reduction
at~$\gp$, and the point $P\in A(L)$ satisfies
\[
  [p](\t-1)P = O.
\]
We would like to conclude that $(\t-1)P=O$, but unfortunately this
need not be true. However, we will show that it is true if we
modify~$P$ by a torsion point.  The following well-known description
of the formal group of an abelian variety at primes of ordinary
reduction provides the crucial information needed to find the required
torsion point.

\begin{theorem}
\label{theorem:ordinaryformalgroup}
Let $K_\gp^\nr$ be the maximal unramifed extension of~$K_\gp$, and
let~$\O_\gp^\nr$ be its ring of integers. Assume that~$A$ has good
ordinary reduction at~$\gp$. Then the formal group~$\Ahat$
of~$A$ is toroidal over~$\O_\gp^\nr$, that is, it is isomorphic
over~$\O_\gp^\nr$ to the formal torus~$\hat\GG_m^g$.
\par
In particular, for all $n\ge1$,
\[
  \Ahat[p^n]
  \cong \hat\GG_m^g[p^n]
  \cong \bfmu_{p^n}^g
\]
as $\Gal(\Kbar_\gp/K_\gp^\nr)$-modules.
\end{theorem}
\begin{proof}
See~\cite[Lemma~3.1]{Baker} or~\cite[Lemma~4.27]{Mazur}.
\end{proof}

Without loss of generality, we may assume that all of the
endomorphisms of~$A$ are defined over~$K$,
%  and that~$K$ contains the CM~field $\End(A)\otimes\QQ$, 
see Proposition~\ref{proposition:reductionstep}(b). Under this assumption,
it is well known (see, for example, \cite[Cor.~2 of Theorem~5]{SerreTate}) 
that the torsion points of~$A$ generate abelian  extensions of~$K$.
To keep the exposition as self-contained as possible, we present a
proof of this fact.

\begin{theorem}
\label{theorem:CMgeneratesabelianexts}
Let~$K$ be a number field, let~$A/K$ be an abelian variety with complex
multiplication. Assume that~$\End_K(A)=\End_\Kbar(A)$.
% and~$\End_K(A)\otimes\QQ\subset K$.
Then $K(A_\tors)\subset K^\ab$, i.e., the torsion points of~$A$ generate
abelian extensions of~$K$.
\end{theorem}

\begin{proof}
Since $\End_K(A)=\End_\Kbar(A)$, it follows from the theory of complex
multiplication that $F = \End_K(A)\otimes\QQ$ is a
CM~field with $[F:\QQ]=2d$, where $d=\dim(A)$.  
If $\ell$ is any prime
number, then one also knows that $V_\ell(A)=T_\ell(A)\otimes\QQ$ has rank one as a module
over $F \otimes \QQ_l$ (see \cite[Lemma 2]{Ruppert}).  Let $G$ be the
image of the natural map $\Gal(\Kbar/K) \to \Aut(V_\ell(A))$.
Since every endomorphism of $A$ is defined over $K$, we get an
injection
\[
G \into (F \otimes \QQ_l)^*,
\]
which implies that $G$ is abelian.  
Since~$\ell$ is an arbitrary prime, the action of~$\Gal(\Kbar/K)$
on~$A_\tors$ is therefore abelian.
\end{proof}

% Alternate proof:
% We observe that~$F$ acts on~$V_\ell(A)=T_\ell(E)\otimes\QQ$, and choosing a basis
% for~$V_\ell(A)$, we can embed~$F$ as a subring of the matrix
% ring~$M_{2d}(\QQ_\ell)$. Similarly, using the action
% of~$\Gal(\Kbar/K)$ on~$T_\ell(A)$, we can map the group
% ring~$\QQ_\ell[\Gal(\Kbar/K)]$ onto a subring~$F'$
% of~$M_{2d}(\QQ_\ell)$. The assumption that~$\End(A)=\End_K(A)$ implies
% that~$F$ and~$F'$ commute. However, the commutant of~$F$ inside
% $M_{2d}(\QQ_\ell)$ is~$F$ itself. Hence~$F'\subset F$, and in
% particular,~$F'$ is commutative. This proves that the action
% of~$\Gal(\Kbar/K)$ on~$T_\ell(A)$ is abelian, and since~$\ell$ is
% an arbitrary prime, the action of~$\Gal(\Kbar/K)$ on~$A_\tors$ is abelian.

We now prepare to state and prove a result (Theorem~\ref{theorem:maintheoremramifiedCMcase}) which, together with 
Theorem~\ref{theorem:maintheoremramifiednonCMcase}, will complete the
proof of our main theorem.
The proof of Theorem~\ref{theorem:maintheoremramifiedCMcase} 
will use the following induction hypothesis.

\begin{center}
  \Hypothesis($f$):\quad
  $\left(\parbox{.75\hsize}{\noindent
    For all fields $K'\subset K^\ab$ whose local conductor satisfies
    \text{$\ord_\gp(\gf_\gp(K'/K))\le f$} and all nontorsion points \text{$P\in A(K')$},
    the height of~$P$ satisfies \text{$\hhat(P) \ge \dfrac{1}{(12p)^{2g}}$}.
   }\right)$
\end{center}

As before, we note that Theorem~\ref{theorem:maintheoremunramifiedcase} shows
that \Hypothesis(0) is true.

\begin{theorem}
\label{theorem:maintheoremramifiedCMcase}
Suppose that $A$ has complex multiplication (over $\Kbar$).
Let $L\subset K^\ab$ be ramified at~$\gp$ and assume that \Hypothesis($f$)
is true for all $f<\ord_\gp(\gf_\gp(L/K))$. Then 
\[
  \hhat(P) \ge \frac{1}{(12p)^{2g}} 
    \qquad\text{for all nontorsion $P\in A(L)$.}
\]
Hence by induction, \Hypothesis($f$) is true for all~$f\ge0$.
\end{theorem}

\begin{proof}
% Recall that the field~$L_\gP$ is a finite abelian extension
% of~$K_\gp=\QQ_p$, so local class field theory tells us that
% $L_\gP=\QQ_p(\z_m)$ for some primitive~$m^{\text{th}}$~root of
% unity. 
Let $m$ be the local conductor $\gf_\gp(L/K)$ of $L/K$ at $\gp$.
Recall that by definition we have $L_\gP\subseteq \QQ_p(\z_m)$, and $m$ is the
smallest positive integer with this property.

Without loss of generality, we may assume that $L$ contains a
primitive $m$th root of unity, since $L(\z_m) \supseteq L$ is again an abelian
extension of $K$ with local conductor $m$.  It follows that $L_\gP = \QQ_p(\z_m)$.

The assumption that~$L/K$ is ramified at~$\gp$ tells us 
that~$p|m$, and indeed~$\t\in\Gal(L/K)$ was chosen to 
generate~$\Gal\bigl(\QQ_p(\z_m)/\QQ_p(\z_{m/p})\bigr)$
(see Remark~\ref{remark:inertiacongruenceremark}).
Write 
\[
  \text{$m=p^km'$ with~$k\ge1$ and~$p\notdivide m'$.}
\]
Then Theorem~\ref{theorem:ordinaryformalgroup} 
and the fact that the formal group has only $p$-power torsion
imply that
\begin{equation}
  \label{equation:torsioninformalgroup}
  \Ahat(L_\gP^\nr)_\tors
  \cong \Ahat[p^k]
  \cong \bfmu_{p^k}^g
\end{equation}
as $\Gal(L_\gP^\nr/K_\gp^\nr)$-modules.  
(We identify $\Ahat(\overline{K}_\gp)$ with the kernel of reduction in
$A(\overline{K}_\gp)$).

We may replace~$L$ by the extension $L(\Ahat[p^k])$.
This is permissible, since from the above isomorphism,
$L(\Ahat[p^k])/L$ is unramified at~$\gp$, and from
Theorem~\ref{theorem:CMgeneratesabelianexts}, the extension~$K(\Ahat[p^k])$
is an abelian extension of~$K$.
\par
Now consider the point~$Q=(\t-1)P$, which by assumption
satisfies~$[p](Q)=O$. The point~$Q$ is in the formal
group~$\Ahat(L_\gP)$, since~$\t$ is in the inertia group. Hence
under the identification~\eqref{equation:torsioninformalgroup},
the point~$Q$ corresponds to a $g$-tuple of $p^{\text{th}}$~roots
of unity,
\[
  (\eta_1,\ldots,\eta_g)\in \bfmu_{p^k}^g.
\]
Also, by construction, the automorphism~$\t$ is a nontrivial element
of the subgroup
$\Gal\bigl(\QQ_p(\z_{p^k})/\QQ_p(\z_{p^{k-1}})\bigr)$, so there is a
\emph{primitive} $p^{\text{th}}$~root of unity~$\xi\in\bfmu_p$ such
that
\[
  \t(\z_{p^k}) = \xi\z_{p^k}.
\]
In other words, the map
\[
  (\t-1):\bfmu_{p^k}\longrightarrow\bfmu_p,
  \qquad
  \z\longrightarrow \t(\z)/\z,
\]
is surjective. Again referring to the
identification~\eqref{equation:torsioninformalgroup}, this shows that
the map
\[
  (\t-1):\Ahat(L)[p^k]\longrightarrow\Ahat(L)[p]
\]
is surjective.
\par
The point $Q=(\t-1)P$ is a point of order~$p$ in~$\Ahat(L)$, 
so this shows that we can find a point~$T\in\Ahat(L)[p^k]$ satisfying
\[
  Q = (\t-1)T.
\]
It follows that $P-T$ is fixed by~$\t$, so $P-T$ is defined over the
fixed field~$L^\t$ of~$\t$. Since~$L^\t$ has strictly smaller
local conductor than does~$L$ at~$\gp$, the induction hypothesis
says that \text{$\hhat(P-T)\ge1/(12p)^{2g}$}. But~$T$ is a
torsion point, so
\[
  \hhat(P-T)=\hhat(P),
\]
which completes the proof of Theorem~\ref{theorem:maintheoremramifiedCMcase},
and with it the proof of our main theorem (Theorem~\ref{theorem:maintheorem}).
\end{proof}

%%%%%%%%%%%%%%%%%%%%%%%%%%%%%%%%%%%%%%%%%%%%%%%%%%%%%%%%%%%%%%%%%%%%%%

\end{document}